\documentclass{amsart}
\usepackage{amsmath}%
\usepackage{amsfonts}%
\usepackage{amssymb}%
\usepackage{graphicx}
\usepackage{psfrag}
\newtheorem{thm}{Theorem}
\newtheorem{prop}{Proposition}
\newtheorem{lem}{Lemma}
\newtheorem*{CL}{$C^r$ Closing Lemma Problem}
\newtheorem*{theoremA}{Theorem \ref{thm:A}}
\newtheorem*{theoremB}{Theorem \ref{thm:B}}
\newtheorem*{theoremC}{Theorem \ref{thm:C}}
\theoremstyle{remark}
\newtheorem*{rem}{Remark}
\begin{document}

\title[On the Closing Lemma problem on the torus]{On the Closing Lemma problem \\ for vector fields of bounded type on the torus}
\author{Simon Lloyd}
\address{School of Mathematics, University of New South Wales, Sydney, Australia}
\email{s.lloyd@unsw.edu.au}

\begin{abstract}
We investigate the open Closing Lemma problem for vector fields on the $2$-dimensional torus. Under the assumption of bounded type rotation number, the $C^r$ Closing Lemma is verified for smooth vector fields that are area-preserving at all saddle points. 
Namely, given such a $C^r$ vector field $X$, $r\geq 4$, with a non-trivially recurrent point $p$, there exists a vector field $Y$ arbitrarily near to $X$ in the $C^r$ topology and obtained from $X$ by a twist perturbation, such that $p$ is a periodic point of $Y$. 

The proof relies on a new result in $1$-dimensional dynamics on the non-existence of semi-wandering intervals of smooth maps of the circle.
\end{abstract}

\maketitle

\section{Introduction}

The $C^r$ Closing Lemma is a longstanding open problem in topological dynamics and was listed by Smale in 2000 as one of the `problems for the next century' \cite{SmaleProblems}. Its importance arises in the fact that a positive solution would lead to many deep results in relation to genericity, stability and bifurcations. The classical $C^r$ Closing Lemma problem is stated as follows:
\begin{CL}
Let $M$ be a smooth compact manifold and let $r\geq 1$ be an integer. Given a $C^r$ vector field $X\in\mathfrak{X}^r(M)$ on $M$, a non-trivially recurrent point $p\in M$ and a neighbourhood $\mathcal{U}\subset \mathfrak{X}^r(M)$ of $X$, does there exist $Y\in \mathcal{U}$ with a periodic orbit passing through $p$?
\end{CL}

The first partial results in this direction are due to Peixoto \cite{Peixoto62} for vector fields on orientable surfaces and to Anosov \cite{Anosov} for uniformly hyperbolic systems. For the $C^1$ topology, this question was answered fully and affirmatively by the celebrated paper of Pugh \cite{PughClosingLemma}. For $r\geq 2$, the problem is open.

There are partial results for $r\geq 2$, for vector fields on the $2$-dimensional torus $\mathbb{T}^2$, the simplest manifold which supports a non-trivial recurrent vector field. 
Gutierrez \cite{GutierrezUnboundedType} proved the $C^r$ Closing Lemma for a large class of vector fields on $\mathbb{T}^2$: the vector fields with `unbounded type' combinatorics. This result was sharpened by Carroll \cite{Carroll92}, and has also been extended to all closed orientable surfaces \cite{Gutierrezhighgenus,AMZ}. Gutierrez and Pires \cite{GutierrezPires} have also demonstrated $C^r$ closing for a class of vector fields on surfaces with an asymptotic contraction property.
In each of these cases, the closing perturbation can take a simple form, a so-called `twist perturbation', supported on an annulus and transverse to the vector field. 

In this article we obtain a positive result for complementary case of `bounded type' vector fields, although all current results for this case are negative. Gutierrez \cite{GutierrezCounterexample} constructs a vector field on $\mathbb{T}^2$, and shows that no local perturbation is a $C^2$ closing perturbation. Carroll showed that, unlike in the unbounded type case, twist perturbations cannot achieve $C^2$ closing for all bounded type vector fields on $\mathbb{T}^2$. The obstruction to closing is that all sufficiently small $C^2$ perturbations produce either another Cantor recurrent set, or else eliminate all one-dimensional recurrent orbits due to the intersection of basins of attractors and repellors. 

We study vector fields on $\mathbb{T}^2$ with singularities in terms of their induced action on a transverse loop $\Sigma$.
The maps of the circle that arise in this way are typically not diffeomorphisms, since if $\Sigma$ intersects the basin $B$ of an attractor (or repellor), then an interval of points in $B\cap\Sigma$ fails to return to the transverse loop in forward (or backward) time. These intervals create plateaus and discontinuities in the induced map.
Thus we instead look at the larger class of order-preserving maps, which have no assumption of continuity, injectivity or surjectivity. 
These maps can possess dynamical features not exhibited by diffeomorphisms, such as `semi-wandering' intervals, whose iterates are disjoint intervals in forward or backward time only.

In Section \ref{semiwandering} we prove the following non-existence of wandering intervals result for order-preserving maps with bounded type rotation number.

\begin{theoremA}
Let $f:S^1\to S^1$ be an order-preserving map of the circle with rotation number of bounded type. 
Suppose that $f$ is $C^1$ except possibly at a finite number of points, and that $\mathop{\mathrm{var}}_{\mathrm{supp}(\mathrm{D}f)} \log \mathrm{D}f$ is finite.

\noindent Then $f$ has either no forward wandering intervals, or else no backward wandering intervals.
\end{theoremA}

In Section \ref{blackcells}, we discuss the features of phase portraits of vector fields on $\mathbb{T}^2$, in particular the grey cells that are characteristic of Denjoy vector fields, and the black cells of Cherry vector fields. We prove a structure theorem for vector fields with bounded type rotation number. 

\begin{theoremB}
Let $X\in\mathfrak{X}^r(\mathbb{T}^2)$ be a bounded type $C^r$ vector field on the torus, $r\geq 4$.
Suppose there are at most finitely many singularities, all hyperbolic, and that the divergence is zero at each saddle point.

\noindent Then $X$ has no grey cells, and any black cells are co-directed.
\end{theoremB}

There are related results on the non-existence of grey cells for vector fields that are repellor-free, and thus automatically have co-directed black cells. Repellor-free vector fields possess backward black cells but no grey cells if the divergence is non-positive at all saddles (see Martens \emph{et al} \cite{MvSdMM}) or else non-negative at all saddles (see Aranson \emph{et al} \cite{Zhuzhoma}). The existence of grey cells in the mixed case, with saddle points with divergence of both signs, is an open problem. 

In Section \ref{closing}, we show that for a vector field with co-directed black cells, it is possible to choose a transverse loop such that the induced map extends to a continuous map of the circle, see Lemma \ref{lem:transversal}. In this case, any $1$-parameter family of $C^r$-twist perturbations contains arbitrarily $C^r$-small closing perturbations.

\begin{theoremC}
Let $X\in\mathfrak{X}^r(\mathbb{T}^2)$, $r\geq 4$, be a bounded type vector field with finitely many singularities all hyperbolic, and suppose the divergence is zero at each saddle. 

\noindent Then for any recurrent point $p$, there exists a loop $\Sigma$ through $p$ and transverse to $X$ such that any $1$-parameter family of $C^r$ twist perturbations along $\Sigma$ contains closing perturbations for $p$ that are arbitrarily $C^r$ near to $X$.
\end{theoremC}

By Proposition \ref{prop:twist}, $C^r$-twist perturbations are sufficient for closing in any vector fields with co-directed black cells, but as Carroll's examples show \cite{Carroll92}, there are vector fields, with contra-directed black cells, for which $C^r$-twist perturbations are ineffective for closing recurrent orbits. 
Carroll's examples use \emph{degenerate} saddle singularities, as does the vector field of Gutierrez with poor $C^2$ closing properties \cite{GutierrezCounterexample}. It is an open problem whether there exist smooth vector fields with contra-directed black cells for which all singularities are hyperbolic. 

\section{Semi-wandering intervals}\label{semiwandering}

Let $S^1$ be the unit circle, with orientation and metric $d$ induced by 
$\pi:\mathbb{R}\to S^1$, the canonical projection $\pi(x)= x\ (\mathrm{mod}\ 1)$. 
We say a function $f:S^1\to S^1$ is a \emph{(strictly) order-preserving map} if there exists a (strictly increasing) 
non-decreasing map $F:\mathbb{R}\to\mathbb{R}$ such that:
\begin{itemize}
\item $F(x+1)=F(x)+1$;
\item $\pi\circ F = f\circ \pi$.
\end{itemize}
Note that order-preserving maps are not assumed to be continuous, injective nor surjective. However, because of their monotonicity, order-preserving maps are continuous except for at most countably many discontinuities of jump type, and are also locally injective, except for countably many intervals of constancy, called \emph{plateaus}.

Denjoy proved that a $C^1$ diffeomorphism with a non-trivially recurrent orbit and with 
derivative of bounded variation is \emph{minimal}: every orbit is dense \cite{Denjoy32}. This result 
is sharp: for arbitrarily small $\epsilon>0$, one can construct a $C^{2-\epsilon}$ 
diffeomorphism which is not minimal due to the presence of a \emph{wandering interval}, 
that is, an interval $I$ whose iterates $f^k(I)$, $k\in\mathbb{Z}$, are pairwise disjoint 
intervals not reduced to a point and the $\omega$-limit set of $I$ is not a periodic orbit. 
In such cases the minimal set is a Cantor set, and is said to be \emph{exceptional}.

Exceptional minimal sets can also occur for order-preserving maps due to the presence 
of ``semi-wandering intervals''. Given a map of the circle $f$, we say an interval 
$I\in\mathbb{T}$ is \emph{forward (backward) wandering} if $f^k(I)$ are disjoint intervals 
not reduced to a point, for integer $k\geq 0$ ($k\leq 0$) and the $\omega$-limit ($\alpha$-limit) 
of $I$ is not a periodic orbit. If an interval $I$ is either forward wandering or 
backward wandering but not both, then it is said to be \emph{semi-wandering}. 
In Theorem \ref{thm:A} we show that for sufficiently smooth maps with bounded type rotation number, 
forward wandering intervals and backward wandering intervals cannot both occur for the 
same map.

Let $R:S^1\to S^1$ denote the rigid rotation $R=R_\alpha:x\mapsto x+\alpha\ (\mathrm{mod}\ 1)$. Recall that for $\alpha$ rational, every orbit is periodic, and that for $\alpha$ irrational, every orbit is dense. In the latter case, information about the ordering of the iterates can be gained from the continued fraction expansion of $\alpha$. Let  $[a_0,a_1,a_2,\ldots]$ denote the continued fraction expansion for $\alpha$, where
$$
\alpha=a_0+\frac{1}{a_1+\frac{1}{a_2+\frac{1}{a_3+\cdots}}}.
$$
Let $(p_n/q_n)$ be the seqence of partial quotients of $\alpha$, obtained by truncating the sequence at the $n$th term. 
The pair of sequences $(p_n)$ and $(q_n)$ generated in this way each obey the same recursive relation
$$
p_n=a_n p_{n-1}+p_{n-2}, \quad q_n=a_n q_{n-1}+q_{n-2}
$$
but have different starting values
$$
p_0=a_0,\ p_1=a_0 a_1+1,\quad q_0=1,\ q_1=a_1.
$$

The sequence $(q_n)$ is called the sequence of \emph{closest return times}, since
$$
d(R^{q_n}(x),x)=\min\{d(R^k(x),x):0<k\leq q_n\},
$$
for any $x\in S^1$. We call the points $R^{q_n}(x)$ the sequence of \emph{closest returns} to $x$.

Note that the sequences $R^{q_{2n}}(x)$ and $R^{q_{2n+1}}(x)$ each converge monotonically to $x$ and from opposite sides. Denote by $I_{2n}$ and $I_{2n+1}$ the decreasing sequences of closed intervals with endpoints $x$ and $R^{2n}(x)$ and $x$ and $R^{2n+1}(x)$ respectively.

The sequence $(a_n)$ from the continued fraction dictates the numbers of ``intermediate returns". Consider the interval $I_n\backslash I_{n+2}$, which has endpoints $R^{q_n}(x)$ and $R^{q_{n+2}}(x)$ from the sequence of closest returns. Any iterates $R^i(x)$, $q_n<i<q_{n+2}$, lying in the interval $I_n\backslash I_{n+2}$ are called \emph{intermediate return points}, and these occur at times $q_n+q_{n+1},q_n+2q_{n+1},\ldots,q_n+(a_n-1) q_{n+1}$: thus they are precisely $a_n-1$ in number.

For the proof of the Theorem \ref{thm:A}, we need to know how much certain collections of intervals overlap.
Given a collection $\mathcal{X}=\{X_i\}_{i\in\mathcal{J}}$ of subsets of a space $X$, the \emph{intersection multiplicity} of $\mathcal{X}$ is defined to be $\textrm{IM}(\mathcal{X})=\max_{x\in X} \#\{i\in\mathcal{J}:x\in X_i\}$.

\begin{lem}\label{lem:Disjointness}
Let $S_n(x)=\{x,R(x),\ldots,R^{q_n-1}(x)\}$, and let $T$ be an interval whose interior is contained in $S^1\backslash S_n(x)$. Then
$$
\mathrm{IM}(\{T,R(T),\ldots,R^{q_n-1}(T)\}) \leq 2 (a_n+1).
$$
\end{lem}
\begin{proof}
The collection of intervals
\begin{eqnarray}\label{TS}
\bigcup_{i=0}^{q_{n}-1} R^i(I_{n-1}) \cup
\bigcup_{i=0}^{q_{n-1}-1} R^i(I_{n})
\end{eqnarray}
cover the circle $\mathbb{T}$, and have pairwise disjoint interiors, see \cite{MS}.

It follows from (\ref{TS}) that the points
\begin{eqnarray}\label{TSptn}
\{x,R(x),\ldots,R^{q_{n-1}+q_{n-2}-1}(x)\}
\end{eqnarray}
partition $S^1$ into intervals of length $|I_{n-1}|$ and $|I_{n-2}|$. Since $q_{n}\geq q_{n-1}+q_{n-2}$, the set $S_n(x)$ is a refinement of partition (\ref{TSptn}), and thus the distance between two adjacent points of $S_{n}(x)$ is also no greater than $|I_{n-2}|$. 

From (\ref{TS}), for any \emph{open} interval $U$ of length $|I_{n-2}|$, the collection $\{ R^i(U) \}_{i=0}^{q_{n-1}-1}$ is pairwise disjoint. 
Thus, allowing for the endpoints, $\textrm{IM}(\{ R^i(\overline{U}) \}_{i=0}^{q_{n-1}-1})\leq 2$.
Since $q_n=a_n q_{n-1}+q_{n-2}< (a_n+1)q_{n-1}$ and since $|T|\leq |I_{n-2}|$, we have the required bound.
\end{proof}

We use an assumption of boundedness of the variation of $\log\mathrm{D}f$ in order to control the distortion of iterates of $f$. 
Recall that for a function $g:U\to \mathbb{R}$, $U\subset S^1$, its \emph{variation} $\mathop{\mathrm{var}}_U(g)$ is defined to be
$$
\mathop{\mathrm{var}}_U(g)=\sup_{x_0<\cdots<x_n}\left\{ \sum_{i=1}^n |g(x_i)-g(x_{i-1})| \right\},
$$
where the supremum is taken over all finite partitions $\{x_i\}_{i=0}^n\subset U$.

\begin{thm}\label{thm:A}
Let $f:S^1\to S^1$ be an order-preserving map of the circle with rotation number of bounded type. 
Suppose that $f$ is $C^1$ except possibly at a finite number of points, and that $\mathop{\mathrm{var}}_{\mathrm{supp}(\mathrm{D}f)} \log \mathrm{D}f$ is finite.

\noindent Then $f$ has either no forward wandering intervals, or else no backward wandering intervals.
\end{thm}
\begin{proof}
We assume, for a contradiction, that exist a forward wandering interval $I_+$ and a backward wandering interval $I_-$.
Let the rotation number of $f$ be $\alpha=[0,a_1,a_2,\ldots]$, with $\liminf a_n =K$ say. Let $\mathcal{N}$ be the set of positive integers such that for all $n\in\mathcal{N}$, $a_n=K$. Note that the set $\mathcal{N}$ is infinite.

Since $\alpha$ is irrational $f$ is semiconjugate to the rigid rotation $R=R_\alpha$: that is, there exists a continuous surjective map $h:S^1\to S^1$ such that $h\circ f(x)=R\circ h(x)$ for all $x\in S^1$. 
Since $R$ is minimal, $h$ maps $I_+$ and $I_-$ to single points. Let $z_+=h(I_+)$ and $z_-=h(I_-)$. We suppose that $I_+$ and $I_-$ are maximal, that is, they are not subsets of larger semi-wandering intervals: otherwise we consider $h^{-1}(z_+)$ and $h^{-1}(z_-)$ instead.
By taking an iterate of $I_+$ if necessary, we may assume that $f^i(I_+)\cap f^{-j}(I_-)=\emptyset$ for every pair of integers $i,j>0$. Since $f$ is $C^1$ except at a finite number of points, we may assume $f^n|_{I_+}$ is $C^1$ for each $n\in\mathbb{N}$, otherwise we replace $I_+$ by some suitably high forward iterate. Similarly, we may assume $f^{-n}|_{I_-}$ is $C^1$ for each $n\in\mathbb{N}$.

Given $n\in\mathcal{N}$, let $q_n$ be the $n$th closest return time for $R$.
Let $\mathcal{F}^+_n=\{f^i(I_+):i=0,\ldots,q_n-1\}$ and $\mathcal{F}^-_n=\{f^{-i}(I_-):i=0,\ldots,q_n-1\}$.
Consider the shortest length interval from the collection $\mathcal{T}_n= \mathcal{F}^+_n \cup \mathcal{F}^-_n$.
Either the shortest length interval in $\mathcal{T}_n$ occurs in the family $\mathcal{F}^+=\{f^i(I_+):i\geq 0\}$ or otherwise in the family $\mathcal{F}^-=\{f^{-i}(I_-):i\geq 0\}$. We denote by $S_+$ the subset of $\mathcal{N}$ non-negative integers for which an interval of $\mathcal{F}_+$ has shortest length among all intervals of $\mathcal{T}_n$, and define $S_-$ correspondingly.
Thus $\mathcal{N}\subset S_+\cup S_-$, and so at least one of $S_+$ and $S_-$ is an infinite set.

Suppose then that $S_+$ is infinite: the case $S_-$ is similar. Fix $n\in S_+$, and suppose $k=k(n)\geq 0$ is such that the shortest length interval of $\mathcal{T}_n$ is the iterate $f^k(I_+)$. Note that by the disjointness of the forward iterates of the forward wandering domain $I_+$, we have that $k(n)\to \infty$ as $n\to\infty$ in $S_+$.

Let $f^{-j_1}(I_-)$ and $f^{-j_2}(I_-)$ be the closest images of $I_-$ on either side of $f^{k}(I_+)$ satisfying $0\leq j_1,j_2 < q_n$. Let $T$ be the shortest closed interval neighbourhood of $I_+$ containing $f^{-j_1-k}(I_-)$ and $f^{-j_2-k}(I_-)$.

By Lemma \ref{lem:Disjointness}, the intersection multiplicity of $T,f(T),\ldots,f^{q_n-1}(T)$ is at most $2a_n+2$, and hence the intersection multiplicity of $T,f(T),\ldots,f^k(T)$ is certainly less than or equal to $2K+2$.

The function $\log \mathrm{D}f$ is real-valued when restricted to the domain $U=\mathrm{supp} (\mathrm{D}f)$. Let $V=\mathop{\mathrm{var}}_U \log \mathrm{D}f$.

Consider $x,y\in \mathop{\mathrm{supp}}\mathrm{D}(f^k)\subset U \cap T$. Then
\begin{eqnarray*}
\log \frac{\mathrm{D}f^k(y)}{\mathrm{D}f^k(x)} & = & \log \frac{\mathrm{D}f(y)\cdot\mathrm{D}f(f(y))\cdots \mathrm{D}f(f^{k-1}(y))}{\mathrm{D}f(x)\cdot\mathrm{D}f(f(x))\cdots \mathrm{D}f(f^{k-1}(x))} \\
& = & \sum_{i=0}^{k-1} \log \mathrm{D}f(f^i(y)) - \log \mathrm{D}f(f^i(x)).
\end{eqnarray*}
Thus, since each point of the circle lies in at most $2(K+1)$ intervals of the collection $\bigcup_{i=0}^k f^i(T)$, we have
$$
\sup_{x,y\in T\cap\mathop{\mathrm{supp}}\mathrm{D}(f^k)} \log \frac{\mathrm{D}f^k(y)}{\mathrm{D}f^k(x)} \leq 2(K+1)V.
$$
Let $\beta:=e^{2(K+1)V}>0$.

By the Mean Value Theorem there exist points
$y_i\in f^{-j_i-k}(I_-)\subset T$, $i=1,2$, and $x\in I_+\subset T$ such that
\begin{eqnarray*}
\mathrm{D}f^{k}(x) & = & \frac{|f^k(I_+)|}{|I_+|}, \\
\mathrm{D}f^{k}(y_i) & = & \frac{|f^{-j_i}(I_-)|}{|f^{-j_i-k}(I_-)|}.
\end{eqnarray*}
Hence, for $i=1,2$,
$$
\beta \geq \frac{|I_+|}{|f^k(I_+)|}\cdot\frac{|f^{j_i}(I_-)|}{|f^{-j_i-k}(I_-)|} \geq \frac{|I_+|}{|f^{-j_i-k}(I_-)|},
$$
where the final inequality follows from our choice of $k$ which ensures $|f^k(I_+)|\geq |f^{-j_i}(I_-)|$.
So for $i=1,2$, we have $|f^{-j_i-k}(I_-)|\geq |I_+|/\beta$.

Hence for any $n\in S_+$ there exists $k(n)$, where $k(n)\to \infty$ as $n\to\infty$, and $m\geq k(n)$, such that
$$
|f^{-m}(I_-)|\geq \frac{|I_+|}{\beta}.
$$
Thus $\liminf_{n\to\infty} |f^{-n}(I_-)|\geq |I_+|/\beta$, which is impossible since $I_-$ is a backward wandering interval and thus has disjoint backward iterates.
\end{proof}

\begin{rem}
In the proof of Theorem \ref{thm:A}, we do not know beforehand whether we will need to use the boundedness of the distortion of high iterates of $f$ or of $f^{-1}$. This bi-directional nature of the argument means that it cannot be generalised in the manner of Yoccoz \cite{Yoccoz} or Martens \emph{et al} \cite{MvSdMM} to the case where $f\sim \pm|x|^\alpha$, $\alpha\geq 1$ at discontinuities and the endpoints of plateaus, and so we must assume that $\log\mathrm{D}f$ is bounded on $\mathrm{supp}(\mathrm{D}f)$.

Thus in order to use Theorem \ref{thm:A} in the proof of Theorem \ref{thm:B}, we must restrict ourselves to vector fields where the divergence at each saddle point is zero.
\end{rem}

\section{Phase portraits of vector fields on the torus}\label{blackcells}

Given a vector field $X\in\mathfrak{X}^r(\mathbb{T}^2)$, $r\geq 2$, we denote by $\Phi^t(x)$ the solution to the ordinary differential equation $\dot{z}=X(z)$ with initial value $z(0)=x\in\mathbb{T}^2$. The map $\Phi:\mathrm{R}\times \mathbb{T}^2\to\mathbb{T}^2$ is called the \emph{flow} associated to $X$. The \emph{forward (or backward) orbit} of a point $p$ is the set $\mathrm{o}^+(p)=\{\Phi^t(p):t\geq 0\}$ (or $\mathrm{o}^-(p)=\{\Phi^t(p):t\leq 0\}$ respectively); the \emph{orbit} of $p$ is the set $\mathrm{o}(p)=\mathrm{o}^+(p)\cup\mathrm{o}^-(p)$. 
A point $p$ (or its orbit $\mathrm{o}(p)$) is \emph{recurrent} if for any neighbourhood $U$ of $p$, $\mathrm{o}(p)\cap U\neq\emptyset$. Simple examples of recurrent orbits are \emph{singularities}, for which $\mathrm{o}(p)=\{p\}$, and \emph{closed orbits}, for which $\mathrm{o}(p)$ is homeomorphic to a circle.
An orbit is said to be \emph{non-trivially recurrent} if it is recurrent but neither a singularity nor a closed orbit. 

If $X\in\mathfrak{X}^r(\mathbb{T}^2)$ has a non-trivially recurrent orbit, and a transverse loop $\Sigma$, we denote by $\rho(X,\Sigma)$ the rotation number of the induced map $\Sigma\to\Sigma$. If we replace $\Sigma$ by a homologically distinct loop $\Sigma'$, then the new rotation number $\rho(X,\Sigma')$ is related to $\rho(X,\Sigma)$ by a linear fractional transformation.
The \emph{rotation orbit} 
$$
\rho(X)=\left\{ \frac{a\rho(X,\Sigma)+b}{c\rho(X,\Sigma)+d}:a,b,c,d\in\mathbb{Z},\quad ad-bc=\pm1 \right\}
$$
is independent of the choice of transverse loop $\Sigma$, and a topological invariant of the vector field $X$, see \cite{Zhuzhoma}. Since a fractional linear transformation acts as a right shift on the coefficients of a continued fraction, the upper and lower limits $\liminf a_n$ and $\limsup a_n$ of the continued fraction of an element $\alpha=[0,a_1,a_2,\ldots]\in\rho(X)$ of the rotation orbit are independent of the choice of $\alpha$. Thus we may describe a vector field as having \emph{(un)bounded type} without reference to any particular transverse loop.

We recall some terminology for describing the orbit structure of vector fields on the torus.
Given a non-trivially recurrent vector field $X\in\mathfrak{X}^r(\mathbb{T}^2)$, the closure of any non-trivially recurrent point is called a \emph{quasiminimal set}. By Maier's estimate \cite{Maier}, the torus can support at most one quasiminimal set: thus if $p$ and $p'$ are non-trivially recurrent points, then $\overline{\mathrm{o}(p')}=\overline{\mathrm{o}(p)} = Q$. 
We say that an invariant set $A\neq Q$ is an \emph{attractor (or repellor)} if there exists a neighbourhood $U\supset A$ such that $\omega(x)=A$ (respectively $\alpha(x)=A$) for all $x\in U$; the maximal such $U$ is called the \emph{basin} of $A$.

A \emph{cell} is a maximal open path-connected invariant set $C$ such that the limit sets $\alpha(x)$ and $\omega(x)$ are independent of $x\in C$. 
A cell $C$ is a \emph{grey cell} if it lies in the complement of $Q$ and satisfies $\alpha(C)=\omega(C)=Q$. Denjoy \cite{Denjoy32} showed that for singularity-free $C^r$ vector fields, grey cells can occur if $r<2$, but if $X$ is $C^1$ with derivative of bounded variation, then $Q=\mathbb{T}^2$ and thus there are no grey cells.

In order to describe the phase portraits of vector fields with singularities, we make use of definitions introduced in the study of Cherry vector fields. 
Cherry \cite{Cherry38} constructed an analytic vector field on the torus with a non-trivially recurrent orbit that has two singularities -- a source and a saddle.
The basin of the source is a cell which lies outside of the quasiminimal set: a cell disjoint from $Q$ that has a limit set other than $Q$ is called a \emph{black cell}.
More specifically, we say a cell $C$ is called a \emph{forward black cell} if $\alpha(C)\neq \omega(C) = Q$, a \emph{backward black cell} if $\omega(C)\neq \alpha(C) = Q$, otherwise $\alpha(C)\neq Q$ and  $\omega(C)\neq Q$ in which case we call $C$ a  \emph{medial black cell}. In each case, the non-quasiminimal limit set, whether an attractor or repellor, is either a node or a limit cycle.

\begin{figure}[htbp]
	\centering
	\psfrag{A}{$C^+$}
	\psfrag{B}{$C^-$}
		\includegraphics[width=0.70\textwidth]{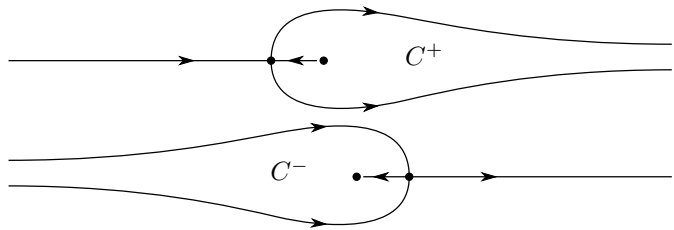}
	\caption{Phase portraits of a forward black cell $C^+$ and a backward black cell $C^-$.}
	\label{fig:BlackCells}
\end{figure}

The phase portrait of a vector field with black cells may be obtained from a singularity-free flow on a torus by `blowing-up' a individual orbits (see \cite{Zhuzhoma}): that is, replacing one or more orbits with either of the phase portraits shown in Figure \ref{fig:BlackCells}.

Cherry constructed a non-trivially recurrent analytic vector field on the torus with a backward black cell, and asked whether such examples can have grey cells. Martens \emph{et al} \cite{MvSdMM} show that repellor-free smooth vector fields with non-positive divergence at all saddles have backward black cells but no grey cells. The techniques of Yoccoz \cite{Yoccoz} show that the same is true if the divergence at all saddles is non-negative, see \cite{Zhuzhoma}. It is an open problem whether grey cells can exist in repellor-free vector fields for which the quasiminimal set contains two saddle points with divergences of opposite signs. 
A vector field is said to have \emph{co-directed} black cells if there are no forward black cells or else no backward black cells; a vector field has \emph{contra-directed} black cells if it has a forward black cell and also a backward black cell.

The case of smooth non-trivially recurrent toral vector fields with both attractors and repellors is also little understood, and no explicit examples have yet been constructed where forward and backward black cells coexist. The following result gives a condition under which the black cells are \emph{co-directed}.

\begin{thm}\label{thm:B}
Let $X\in\mathfrak{X}^r(\mathbb{T}^2)$ be a bounded type $C^r$ vector field on the torus, $r\geq 4$.
Suppose there are at most finitely many singularities, all hyperbolic, and that the divergence is zero at each saddle point.

\noindent Then $X$ has no grey cells, and any black cells are co-directed.
\end{thm}
\begin{proof}
We assume $X$ has singularities, otherwise there are no black cells, and by Denjoy's theorem \cite{Denjoy32} there are no grey cells.

Since $X$ has a non-trivially recurrent orbit, we may construct a $C^r$ loop $\Sigma$ transverse to $X$, see \cite{PalisDeMelo}. 
Let $f:\Sigma\to\Sigma$ be an order-preserving extension to $\Sigma$ of the induced map. As $X$ is $C^4$ and the divergence at each saddle points is zero, it follows that $\mathop{\mathrm{var}}_{\mathrm{supp}(\mathrm{D}f)}\log \mathrm{D}f$ is finite, see for example \cite{MvSdMM}.

Assume for a contradiction that there exist a forward black cell $C^+$ and a backward black cell $C^-$.  
Then $C^+\cap\Sigma$ contains an interval $I^+$ that is a forward wandering interval for $f$, and $C^-\cap\Sigma$ contains a backward wandering interval $I^-$ for $f$. As this contradicts Theorem \ref{thm:A} applied to $f$, it follows that $X$ has either no forward black cells or else no backward black cells.

Assume now that $X$ has a grey cell $C$. Then $C\cap\Sigma$ contains backward and forward wandering intervals $I^-$ and $I^+=f(I^-)$ for $f$. We again obtain a contradiction to Theorem \ref{thm:A}, and so $X$ has no grey cells.
\end{proof}

\section{Closing Lemma results}\label{closing}

When studying vector fields on the torus in terms of the return map to a transverse loop, it is common to consider `twist perturbations'.
Let $U$ be an annular neighbourhood of a transverse loop $\Sigma$ that is free of singularities. A \emph{family of $C^r$ twist perturbations} is a continuous family of $C^r$ vector fields $Y:[0,a_0]\to \mathfrak{X}^r(\mathbb{T}^2)$ such that $Y_0\equiv X$ and for $a>0$, $\mathrm{supp}(Y_a-X)= U$ and $Y_a|_{\mathrm{int}\ U}$ is transverse to $X|_{\mathrm{int}\ U}$.

We say $X$ admits \emph{$C^r$ closing by twist perturbations along $\Sigma$} if for any non-trivially recurrent point $p\in\Sigma$, any annular neighbourhood $U$ of $\Sigma$ and any family $\{Y_a\}_{0\leq a<a_0}$ of $C^r$ twist perturbations, there exists $(a_n)\downarrow 0$ and points $x_n\in\Sigma$, $(x_n)\to p$, such that $x_n$ lies on a closed orbit of $Y_{a_n}$. 
As noted by Pugh \cite{PughClosingLemma}, in such a case we obtain an affirmative answer to the $C^r$ Closing Lemma problem as follows. Let $Z_n:\mathbb{T}^2\to \mathbb{T}^2$ denote the translation of the torus which takes $p$ to $x_n$.
If $X$ admits $C^r$ closing by twist perturbations along a transverse loop $\Sigma$, then by setting $X_n(x):=Y_{a_n}(Z_n(x))$, we obtain a sequence of vector fields with closed orbits through $p$ and satisfying $\|X-X_n\|_{C^r}\to 0$.

The return map induced by a vector field on a transverse loop $\Sigma$ typically exhibits plateaus, caused points lying in the basin of an attractor whose forward orbits do not return to $\Sigma$, and also discontinuities caused by points in basins of repellors whose backward orbits do not return to $\Sigma$. We begin by proving that when black cells are co-directed, it is possible to find a transverse loop such that the induced map lacks one or the other of these features. 

\begin{lem}\label{lem:transversal}
Let $X\in\mathfrak{X}^r(\mathbb{T}^2)$, $r\geq 1$, be a vector field with finitely many singularities all hyperbolic, and suppose $X$ has a non-trivially recurrent point $p$. If there are no forward black cells, then there exists a transverse $C^r$ loop $\Sigma$ through $p$ such that the induced map $P:\Sigma\to\Sigma$ extends to a continuous map.
\end{lem}
\begin{proof}
It is enough to show that we can construct a tranverse loop disjoint from all repellor basins.
Since there are no backward black cells, each repellor basin consists of a disjoint union of medial black cells.
Let $Q$ be the quasiminimal set and let $H$ denote the union of the medial black cells. 
Since $\overline{H}\cap Q$
consists only of saddle points and separatrices, we have that $\mathrm{o}(p)$ is disjoint from $\overline{H}$. Since $p$ lies in the open set $\mathbb{T}^2\backslash \overline{H}$, we may take a small open disc $U$ around $p$ disjoint from $\overline{H}$. Let $p_0,p_1$ be the first exit and first re-entry points of $\mathrm{o}^+(p)$ to $U$. Let $\gamma=\mathrm{o}^+(p_0)\cap\mathrm{o}^-(p_1)$ be the segment of orbit from $p_0$ to $p_1$, and let $L$ be an open tubular neighbourhood of $\gamma$ disjoint from $\overline{H}$. For $\epsilon$ sufficiently small, the orbit of $p$ under the vector field $X+\epsilon X^\perp$ exits $U$ at a point $q_0$ and re-enters $U$ at a point $q_1$ with $\ell=\mathrm{o}^+(p_0)\cap\mathrm{o}^-(p_1)$ contained in $L$. By connecting $q_1$ to $q_1$ with a $C^r$ path $\ell'\subset U$ through $p$ transverse to $X$, we obtain a $C^r$ loop $\Sigma=\ell\cup \ell'$ that is transverse to $X$ and disjoint from $\overline{H}$.

Let $P:\Sigma\to\Sigma$ be the induced return map. By construction, $\Sigma$ does not intersect the basins of any repellors. Thus $P$ extends to a continuous map $\Sigma\to\Sigma$.
\end{proof}

\begin{lem}\label{lem:periodicpoints}
Let $f_0:S^1\to S^1$ be a continuous map of the circle that is strictly increasing except for a finite number of plateaus, and let $p\in S^1$ be a non-trivially recurrent point whose orbit accumulates on $p$ from both sides. 
Suppose $f_{\epsilon}:=h_{\epsilon}\circ f_0$, where $h_\epsilon:S^1\to S^1$ is a continuous family of homeomorphisms, with lifts $H_\epsilon:\mathbb{R}\to \mathbb{R}$ such that $H_0=\mathrm{Id}$ and $\epsilon\mapsto H_\epsilon(x)$ is continuous, non-decreasing at each $x$ and strictly increasing at $p$.

Then there exists a sequence $s_n\downarrow 0$ such that for each $n\in\mathbb{N}$, $f_{s_n}$ has a periodic point $x_n\to p$ and moreover the orbit of $x_n$ is disjoint from the plateaus of $f_{s_n}$.
\end{lem}
\begin{proof}
Since $p$ is non-trivially recurrent, $\rho(f_0)$ is irrational. Moreover, since $\epsilon\mapsto H_\epsilon(p)$ is strictly increasing, the mapping
$\epsilon\mapsto \rho(f_\epsilon)$ is strictly increasing at $\epsilon=0$. As $f_0$ is continuous, the mapping $\epsilon\mapsto \rho(f_\epsilon)$ is continuous, and thus we can find a sequence $(s_n)\downarrow 0$ such that $\rho(f_{s_n})=p_n/q_n$.

If for arbitrarily large $n$, there is a periodic orbit through $p$ disjoint from all plateaus of $f_{s_n}$ we are done. 
Thus, shifting $(s_n)$ if necessary,  we suppose this is not the case. Let $I_n$ denote the shortest interval containing $p$ and two distinct periodic points $x_n,y_n$ (of period $q_n$) of $f_{s_n}$ as its endpoints. Since $\epsilon\mapsto \rho(f_\epsilon)$ is locally constant at rational values, we may perturb each $s_n$ slightly, leaving $\rho(f_{s_n})$ unchanged, and ensure that the orbits of $x_n$ and $y_n$ do not intersect the endpointa of any plateau of $f_{s_n}$. It then follows from the Intermediate Value Theorem applied to $(f_{s_n})^{q_n}$ that at least one of the periodic endpoints, $x_n$ say, has an orbit disjoint from all plateaus of $f_{s_n}$.

It remains to show that $x_n\to p$. Assume for a contradiction that there exists an interval $J\ni p$, $J\subset I_n$ for arbitrarily large $n$. Since $I_n$ has period $q_n\to\infty$, it follows that $\{(f_{s_n})^k(J):0\leq k<q_n\}$ are pairwise disjoint intervals. Thus $p$ is contained in a wandering interval for $f_0$, which contradicts the assumption that $\mathrm{o}^+(p)$ accumulates on $p$ from both sides.
\end{proof}

\begin{prop}\label{prop:twist}
Let $X\in\mathfrak{X}^r(\mathbb{T}^2)$, $r\geq 1$, be a vector field with co-directed black cells and a non-trivially recurrent point $p$. 
Then there exists a loop $\Sigma$ along which $X$ admits $C^r$ closing by twist perturbations.
\end{prop}
\begin{proof}
Since $X$ is has co-directed black cells, by Lemma \ref{lem:transversal} there exists a transverse loop $\Sigma$ through $p$ such that the induced map to $\Sigma$ extends to a continuous map $P:\Sigma\to\Sigma$. 
Let $U$ be a singularity-free annular neighbourhood of $\Sigma$ and let $\{Y_a\}_{0\leq a<a_0}$ be a family of $C^r$ twist perturbations.
Reducing $a_0$ if necessary, we may assume $\Sigma$ is transverse to $Y_a$ for all $a\in [0,a_0)$.

Let $f_a:\Sigma\to\Sigma$ be the map induced by $Y_a$.
Then the family $\{f_a\}_{0\leq a<a_0}$ satisfies the hypotheses of Lemma \ref{lem:periodicpoints}, and thus there exist sequences $x_n\to p$ and $(a_n)\downarrow 0$, where $x_n\in\Sigma$ is a periodic point of $f_{a_n}$. Moreover, since the orbit of $x_n$ does not intersect any plateaus of $f_{a_n}$, we have that the orbit of $x_n$ under $Y_{a_n}$ is a closed orbit.
\end{proof}

\begin{thm}\label{thm:C}
Let $X\in\mathfrak{X}^r(\mathbb{T}^2)$, $r\geq 4$, be a bounded type vector field with finitely many singularities all hyperbolic, and suppose the divergence is zero at each saddle. 

\noindent Then for any recurrent point $p$, there exists a loop $\Sigma$ through $p$ and transverse to $X$ such that any $1$-parameter family of $C^r$ twist perturbations along $\Sigma$ contains closing perturbations for $p$ that are arbitrarily $C^r$ near to $X$.
\end{thm}
\begin{proof}
The result follows directly from Theorem \ref{thm:B} and Proposition \ref{prop:twist}.
\end{proof}

\begin{rem}
In Theorems \ref{thm:B} and \ref{thm:C} we assume $X$ is $C^4$ because of the resonance of the eigenvalues of $X$ at the 
saddle points (see for example \cite{Stowe,MvSdMM}).
Both theorems remain true when $X$ is $C^2$, so long as there are $C^2$ linearising coordinates for each saddle point.
\end{rem}

\section*{Acknowledgements}
The author wishes to thank S. van Strien for discussions on this topic.

\end{document}